\documentclass{article}      

\usepackage{amssymb}
\usepackage{amsmath}
\usepackage{latexsym}
\usepackage[mathscr]{euscript}
\usepackage{graphicx}
\usepackage{amscd}
\usepackage{amsthm} 
\usepackage{enumerate}
\usepackage{wrapfig}
\newtheorem{theorem}{Theorem}
\newtheorem{corollary}[theorem]{Corollary}
  
\newtheorem{proposition}[theorem]{Proposition}
\newtheorem{conjecture}[theorem]{Conjecture}

\newcommand{\bc}{\mathbb{C}}
\newcommand{\bp}{\mathbb{ P}}
\newcommand{\bz}{\mathbb{Z}}
\newcommand{\br}{\mathbb{R}}
\newcommand{\bq}{\mathbb{Q}}

\newcommand{\p}{\partial}

\newcommand{\cs}{\mathcal{S}}
\newcommand{\cg}{\mathcal{G}}
\newcommand{\cb}{\mathcal{B}}
\newcommand{\cp}{\mathcal{P}}
\newcommand{\ch}{\mathcal{H}}
\newcommand{\cm}{\mathcal{M}}

\newcommand{\ct}{\mathcal{T}}

\newcommand{\ca}{\mathcal{A}}

\newcommand{\hk}{\hookrightarrow}
\newcommand{\bg}{\bigskip}
\newcommand{\med}{\medskip}
\newcommand{\la}{\longrightarrow}
\newcommand{\bfl}{\begin{flushleft}}
\newcommand{\efl}{\end{flushleft}}

\newcommand{\eps}{\epsilon}

\newcommand{\bcp}{\bc \bp}

 \newcommand{\cy}{\mathcal{Y}}

\newcommand{\xr}{\xrightarrow}

\newcommand{\blt}{\bullet}

\newcommand{\casd}{\ca_{sd}}

\newcommand{\bgg}{B\Gamma_{\infty, 1}}

 \DeclareMathOperator{\Diff}{Diff}

\newcommand{\dgn}{\Diff^+(\Sigma_{g,n}, \p \Sigma_{g,n})}
 
 \newcommand{\dsg}{\Diff^+(\Sigma_g)}

\begin{document}

 \title{Stability phenomena in the topology of moduli spaces}

 \author {
Ralph L. Cohen  \thanks{The   author was  partially supported by  NSF grant  DMS-0603713}
   \\ Dept. of Mathematics\\ Stanford University\\  
 Stanford, CA 94305
  }
 \date{\today}
\maketitle  
\begin{abstract} The recent proof by Madsen and Weiss of Mumford's conjecture on the stable cohomology of moduli spaces of Riemann surfaces,
was a dramatic  example of an important stability theorem about the topology of moduli spaces.  In this article we give a survey of families of classifying spaces and moduli spaces where 
  ``stability  phenomena" occur in their topologies.   Such  stability theorems have been proved in many situations in the history of topology and geometry, and the payoff  has often been  quite remarkable.  In this paper we discuss classical stability theorems such as the Freudenthal suspension theorem, Bott periodicity, and Whitney's embedding theorems.  We then discuss more modern examples such as those involving  configuration spaces of points in manifolds,  holomorphic curves  in complex  manifolds,  gauge theoretic moduli spaces,  the stable topology
of general linear groups,  and pseudoisotopies of manifolds.  We then discuss the stability theorems regarding the moduli spaces of Riemann surfaces:  Harer's stability theorem on the cohomology of moduli space, and the Madsen-Weiss theorem, which proves a generalization of Mumford's conjecture.  We also
 describe Galatius's recent theorem on the stable cohomology of automorphisms of free groups.    We end by speculating on the existence of general conditions in which one might expect these stability phenomena to occur.   \end{abstract}

 \tableofcontents

 \section*{Introduction}    In the last sixty years, the notions  of classifying space  and  moduli space  have played   central roles in the development of  topology and geometry.
 These are spaces that encode the basic topological or geometric structure to be studied, and therefore the topology of these spaces  naturally have been a subject of intense interest. 
 Probably the most fundamental among them are the moduli spaces of Riemann surfaces of genus $g$, $\cm_g$.  In a dramatic application of algebraic topological methods to algebraic geometry, Madsen and Weiss recently proved a well known conjecture of Mumford regarding the \sl stable \rm cohomology of moduli space  \cite{madsenweiss}.
 Namely,  Mumford described a ring homomorphism  from a graded  polynomial algebra over the rationals, to the cohomology of moduli space with rational coefficients,
 $$
 \bq [\kappa_1, \kappa_2, \cdots \kappa_i, \cdots]  \la H^*(\cm_g;  \bq),
 $$
 and conjectured that it is an isomorphism when the genus $g$ is large with respect to the cohomological grading.
 Here $\kappa_i$ is the Miller-Morita-Mumford  canonical class,  and has grading $2i$.   In \cite{madsenweiss} Madsen and Weiss described a homotopy
 theoretic model for the \sl stable \rm moduli space, $\cm_\infty$, and in so doing, not only proved Mumford's conjecture, but also gave an implicit model for the stable cohomology of moduli space with \sl any \rm coefficients.  Using this explicit model, Galatius \cite{galatius1} calculated this stable cohomology explicitly, when the coefficients are $\bz/p$ for $p$ any prime, and in so doing uncovered a vast amount of previously undetected torsion in the stable cohomology of moduli space.

 \med
 The Madsen-Weiss theorem can be viewed as one of the most recent  examples  of a stability theorem  regarding the topology of classifying spaces or moduli spaces.  The purpose of this paper is to give a survey of these types of theorems and their applications to a broad range of topics in topology and geometry.

 \med
 Stability theorems are results regarding \sl families \rm of classifying spaces or moduli spaces.  These spaces are typically indexed by
   some geometrically defined quantity, such as the degree of a map,  the rank of a bundle,  the genus of a curve, or a characteristic number. We refer to this number as the ``degree".  (In the case of the moduli spaces of curves, this indexing degree is the genus of the curve.)  We let $\cm_d$ be the moduli space corresponding to degree $d$. 
   
   Two basic questions about the topology of these spaces  naturally occur, and seeing how they are addressed in a variety of examples is the basic theme of this paper.
   
 \begin{itemize}
 \item \sl{Stability Question 1.} \rm
   How does the topology of the moduli spaces change as the degree changes?  Is there a ``stability range" for their homology or homotopy groups?  By this we mean a function $r(d)$ which is  an unbounded and nondecreasing function of the degree $d$,  with the property that the $k^{th}$ homology and/or homotopy group of   $\cm_d$ and  $\cm_{d+1}$ are isomorphic so long as $k < r(d)$.  
   \item \sl{Stability Question 2.}  \rm
     Is there a naturally defined, more easily accessible limiting homotopy type, as the degree gets large? If so, calculate this ``stable homotopy type" as explicitly as possible.
     \end{itemize}
 
 In this survey article we discuss a variety of examples of families of classifying spaces and moduli spaces where these questions have been addressed.
  Different techniques have been
 used to study these questions, but as we hope to point out,  there are common themes among these techniques.

\med
We organize this survey in the  following way.  In section one, we discuss classical stability theorems, including the Freudenthal suspension theorem,  Bott's periodicity theorem, and Whitney's embedding theorems.   In sections 2 through 4 we discuss more modern stability theorems, including those dealing with configuration spaces of points in manifolds,  holomorphic curves  in complex  manifolds,  gauge theoretic moduli spaces,  the stable topology
of general linear groups,  and pseudoisotopies of manifolds.   In section 5 we discuss the background of Mumford's conjecture, including the stability theorem for mapping class groups of Harer.  We then discuss the Madsen-Weiss theorem in some detail, and also describe similar  theorems regarding automorphism groups of free groups. This  includes stability theorems of Hatcher and  Vogtmann,   and the recent theorem of Galatius about the stable cohomology of automorphisms of free groups.  We end with  a discussion in  section 6 regarding potential research questions whose goal is to find general criteria under which stability theorems hold (and do not hold).

   \section{Classical stability theorems}
   
   \subsection{The Freudenthal suspension theorem}
   Probably the oldest example of a stability theorem in topology and geometry, proved in 1938,  is the ``Freudenthal suspension theorem" \cite{freudenthal}.  Let $\Omega^dS^d$ be the space of self maps of the sphere $S^d = \br^d \cup \infty$ that fix the basepoint at infinity.  By the adjoint construction, there is a natural identification of homotopy groups,
   $$
   \pi_q \Omega^dS^d\cong \pi_{q+d}S^d.
   $$
   Moreover there  is a natural  ``suspension" map
   $$
   \Sigma : \Omega^dS^d \to \Omega^{d+1}S^{d+1}
   $$
   defined as follows.   Let $X$ be any space with a  fixed basepoint $x_0 \in X$.  The suspension of $X$, written $\Sigma X$ is the quotient
   $$
   \Sigma X = S^1 \times X /( \infty \times X) \cup (S^1 \times x_0).
   $$
   This construction is natural, in the sense that if one has a (basepoint preserving) map $f: X \to Y$, then one has an induced ``suspension map",
   $\Sigma f : \Sigma X \to \Sigma Y$ defined by
    $\Sigma f (t,x) = (t, f(x)).$
    
    There is a natural identification (homeomorphism) of $\Sigma S^d \cong S^{d+1}$,  with respect to which  the suspension construction defines the map $\Sigma : \Omega^dS^d \to \Omega^{d+1}S^{d+1}$. 
 The following is Freudenthal's basic theorem:
 
 \begin{theorem}\label{freu1}
 The suspension map $\Sigma : \Omega^dS^d \to \Omega^{d+1}S^{d+1}$ induces an isomorphism in homotopy groups 
 $$
\Sigma_*:  \pi_q(\Omega^dS^d) \xr{\cong} \pi_q(\Omega^{d+1}S^{d+1}) 
 $$
 for $q < d-1$.  It is a surjection for $q = d-1$.  In other words $\Sigma_*: \pi_r S^d \to \pi_{r+1}S^{d+1}$ is an isomorphism
 for $r \leq 2(d-1)$, and is a surjection for $r = 2d-1$.  
 \end{theorem}
 
 \med
 Notice that this result can be viewed as answering Stability Question 1  in this setting.  This theorem has the following generalization.   Let $X$ be any $k$-connected
 space with a distinguished basepoint, $x_0 \in X$.  That is,  $\pi_r X = 0$ for $r \leq k$.  Let $\Omega^d X$ denote the space of continuous maps $\alpha : S^d \to X$ that take the basepoint $\infty \in S^d$ to $x_0$.     Suspending defines a map
 $$
 \Omega^d X \to \Omega^{d+1}\Sigma X.
 $$
 The following gives a generalization of the above theorem:
 \begin{theorem} \label{freu2}
 $$ 
\Sigma_*:  \pi_q ( \Omega^d X ) \to \pi_q(\Omega^{d+1}\Sigma X)
$$
is an isomorphism for $q \leq 2k -d$, and is surjective for $q = 2k-d+1$.  
In other words,
$$
\Sigma_* : \pi_j X \to \pi_{j+1}(\Sigma X)
$$
is an isomorphism for $j \leq 2k$ and is a surjection for $j = 2k+1$. 
\end{theorem}

 As mentioned, these results can be viewed as answers to Stability Question 1 in this context.  To address Stability Question 2, one considers
 the limiting space, $Q(X) = \lim_{n\to \infty} \Omega^n \Sigma^n X$. The homotopy groups of $Q(X)$, are the stable homotopy groups of $X$,
 $$
 \pi_q(Q(X)) = \lim_{n\to \infty} \pi_{q+n}\Sigma^n (X)  = \pi_q^s(X).
 $$
 While these stable homotopy groups are notoriously difficult to compute, they do have a significant advantage over  the unstable homotopy groups.  Namely, the functor $X \to \pi_*^s(X)$ is a (reduced) generalized homology theory,  in that it satisfies the Eilenberg-Steenrod axioms.  In particular the excision axiom holds for the stable theory, but does not hold for unstable homotopy groups.   Over the years this has allowed for a variety of powerful calculational techniques. An important one, for example, is the spectral sequence of Atiyah and Hirzebruch that approximates $\pi_*^s (X)$ by the homology groups, $H_*(X; \pi_*^s)$, where the coefficients,  $\pi_*^s$ are the stable homotopy groups of spheres.

 \subsection{Whitney's Embedding Theorem}
 The classical Whitney Embedding Theorem can be viewed as a stability theorem for the moduli space of smooth submanifolds  of $\br^\infty$ of a given diffeomorphism type.
 More specifically,  let
  $M^n$ be a closed $n$-dimensional smooth manifold.  Let  $Emb(M^n, \br^N)$   be  the space of smooth embeddings $e : M^n \hk \br^N$.  This embedding space is topologized using the compact-open topology.  Whitney's basic embedding theorem \cite{whitney}  is the following.
  
  \begin{theorem}\label{whitney} For $N \geq 2n$,   $Emb(M^n, \br^N)$ is nonempty.   For $N \geq 2(n+k)$,  the space $Emb(M^n, \br^N)$ is $(k-1)$-connected.  That is,
  the homotopy groups,
  $$
  \pi_i(Emb(M^n, \br^N)) = 0
  $$
  for $i \leq k-1$.
  \end{theorem}
  
  Notice that in the above theorem for $k = 1$,  it says that $Emb(M^n, \br^{2n+2})$ is connected; i.e. any two embeddings are isotopic.  The fact that  $ \pi_1(Emb(M^n, \br^{2n+4})) = 0 $  can be interpreted to say that not only are any two embeddings isotopic, but any two isotopies can be deformed to each other by a one-parameter family of isotopies.   Taking the limit as $N \to \infty$, one has that $Emb(M^n, \br^\infty)$ is weakly contractible (i.e. all of its homotopy groups are zero). Indeed it   can easily be shown that this space is contractible, which can be interpreted as saying that not only are any two
  embeddings isotopic, but that there is a contractible family of choices of isotopies between them.
  
  \med
  The diffeomorphism group $\Diff (M)$ acts freely on the embedding spaces, $Emb(M^n, \br^N)$.  The action also is known to admit  slices, which implies that the projection
  onto the quotient, which we call  $\cm_N (M),$ is a fiber bundle.  We can think of $\cm_N (M)$ as the moduli
  space of submanifolds of $\br^N$ that are diffeomorphic to $M$.  As a consequence of Whitney's theorem, one has the answer to Stability Question 1 in this context.
  
  \begin{corollary}\label{whit-stab1}  The linear inclusion  $\br^N \hk \br^{N+1}$ induces a ``gluing map"
  $$
  \cm_N(M^n) \to \cm_{N+1}(M^n)
  $$
  which induces an   isomorphism  in homotopy groups in dimensions less than $\left[\frac{N-2n}{2}\right]$, and is surjective in dimension $\left[\frac{N-2n}{2}\right]$.
  \end{corollary}

  \med
  Now by letting $N \to \infty$, Whitney's theorem also supplies an answer to Stability Question 2 in this setting.  Namely,  since Whitney's theorem implies that  the total  space of the bundle $Emb(M^n, \br^\infty) \to Emb(M^n, \br^\infty)/\Diff (M^n)  = \cm_\infty (M) $ is weakly contractible,   the moduli space can be taken to be the classifying space of the diffeomorphism group,
  $$
 \cm_\infty (M) \simeq B\Diff (M^n).
  $$
 
 This observation can be interpreted in the following way. Consider the moduli space with one marked point,
 $$
 \cm_{\infty, 1} (M) = \{(N, x), \, \text{where} \, N \subset \br^\infty  \, \text{is diffeomorphic to M}, \, \text{and} \, x \in N.\}
 $$
 The projection map
\begin{align}
 p :  \cm_{\infty, 1} (M) &\to \cm_\infty (M) \notag \\
 (N, x) & \to N \notag
 \end{align}
 is a fiber bundle whose fiber is $M$.  It is referred to as the ``canonical" $M$-bundle over $\cm_\infty (M)$.  The following interpretation of $\cm_\infty (M)$ as the classifying space $B\Diff (M)$ has been   used in an important way by Madsen and Weiss in their proof of Mumford's conjecture on the stable homology of the moduli space of curves \cite{madsenweiss}, as well as in the study of cobordism categories \cite{GMTW}.   
  
  \begin{proposition}\label{whit-stab2}  The  stable  moduli space of manifolds diffeomorphic to $M$, $\cm_\infty (M^n)$,  classifies fiber bundles with fiber $M^n$.
  That is, for a space $X$ of the homotopy type of a $CW$-complex, there is a bijective correspondence,
  $$
 \phi :  [X, \cm_\infty (M)] \xr{\cong}     Bdl_M(X)
  $$
  where the left hand side is the set of homotopy classes of maps, and the right hand side is the set of  isomorphism classes of fiber bundles over $X$ with fiber $M$
   and structure group $\Diff (M)$.  The correspondence $\phi$ assigns to a map $f : X \to \cm_\infty (M)$ the pull-back of the canonical bundle, $f^*(\cm_{\infty_,1} (M))$.
  \end{proposition}
  
  \subsection{Bott periodicity}
  
  In \cite{bott} R. Bott proved his famous periodicity theorem on the ``stable" homotopy type of Lie groups.  Primarily this is a theorem about the homotopy type of the orthogonal groups  and unitary groups $O(n)$ and $U(n)$ as $n$ gets large.   These results can be interpreted as stability results about the moduli space of vector spaces, in the following way.
    
    Let $Gr_k(\bc^n)$ be the Grassmannian of $k$-dimensional complex subspaces of $\bc^n$.    By increasing $n$, one can consider the infinite Grassmannian
    $Gr_k(\bc^\infty)$.  In analogy with the above discussion about embeddings of manifolds, this Grassmannian can be viewed as the quotient,
    $$
    Gr_k (\bc^\infty) = Mono (\bc^k, \bc^\infty)/U(k)
    $$
    where $Mono (\bc^k, \bc^\infty)$ is the space of linear monomorphisms that preserve the Hermitian inner product (the ``Stiefel manifold").    This space is acted upon freely by 
    the unitary group, $U(n)$, and the quotient space,  $Gr_k (\bc^\infty)$ can be viewed as the ``moduli space" of $k$-dimensional complex vector spaces.  
   Since $Mono (\bc^k, \bc^\infty)$ is contractible, this space is a model for the classifying space, $BU(k)$, which classifies $k$-dimensional complex vector bundles.
   (See \cite{milnorstasheff} for a thorough discussion.)

  Given a $k$-dimensional space $V \subset \bc^\infty$, then crossing with a line gives $V \times \bc \subset \bc^\infty \times \bc$.  Choosing
  a fixed  isomorphism $\bc^\infty \times \bc \cong \bc^\infty$   defines a ``gluing" map,
  $$
g_k : Gr_k(\bc^\infty) \to   Gr_{k+1}(\bc^\infty).
$$
It is well known (see \cite{milnorstasheff}) that this map is homotopy equivalent to the unit sphere bundle
$$
\pi_{k+1} : S(\gamma_{k+1}) \to Gr_{k+1}(\bc^\infty)
$$
where $S(\gamma_{k+1}) = \{ (W, w) \, : \, W \in Gr_{k+1}(\bc^\infty), \, \text{and} \, w \in W \, \text{with} \, |w|=1.\}$  Since the fiber of $\pi_{k+1}$ is the sphere $S^k$, one has the following answer to Stability Question 1 in this context:

\begin{proposition}  The gluing map
 $$
g_k :BU(k) \to   BU(k+1) 
$$
induces an isomorphism in homotopy groups in dimensions less than $k$,and is a surjection in dimension $k$.
\end{proposition}

In this setting, Bott's theorem, one of the most important theorems in topology in the twentieth century,  can be viewed as an answer to Stability Question 2.   Let $BU = \lim_{k\to \infty} BU(k)$ be the (homotopy) colimit of the gluing maps
$g_k$.

\begin{theorem} (Bott periodicity \cite{bott}) \label{bott}
$$
\pi_q(\bz \times BU)  \cong \begin{cases} \bz  \quad \text{if} \quad q \, \text{is even} \\
0    \quad \text{if} \quad q \, \text{is odd}  \end{cases}
$$ In particular  these homotopy groups are periodic, with period 2.
 
 If $BO$ is defined similarly (using real Grassmannians), then
 $\pi_q(\bz \times BO)$ is periodic of period 8,  and  the first eight homotopy groups (starting with dimension 0) are given by   $\bz, \, \bz/2, \, \bz/2, \, 0, \, \bz, 0, \, 0, \, 0.$
\end{theorem}

\section{Configuration spaces, permutations, and braids}

\subsection{Configurations of points in a manifold}
Let $M$ be a manifold, and $F_k(M) \subset M^k$ be the space of $k$-distinct points in $M$.  There is a natural free action of the symmetric group, $\Sigma_k$,
and we let $C_k(M)$ be the orbit space,
$$  C_k(M) = F_k(M)/\Sigma_k.
$$
$C_k(M)$ is the moduli space of $k$ - points  (or particles) in $M$ and has proven  extremely important in a variety of applications in topology, geometry, and physics.  Generalizing results of Segal \cite{segal},  McDuff  proved  results that answer Stability Questions 1 and 2 in this context.   
These results can be described as follows.   Let $M$ be a smooth, connected, open $n$-dimensional manifold that is the interior of a compact manifold with boundary, $\bar M$.  
Let $p: TM \to M$ be the tangent bundle, and let $T^\infty M \to M$ be the associated  $S^n$-bundle obtained by taking the fiberwise one point compactification of $TM$.  So the fiber of $p : T^\infty M \to M$ at $x \in M$, is the compactified tangent space, $T_xM \cup \infty$.      Let $\Gamma M$ be the space of  smooth sections of $T^\infty M$  that have compact support.  Such a section has a degree.  We write $\Gamma (M) = \coprod_{k\in \bz} \Gamma_k (M)$, where $\Gamma_k (M)$ are the sections of degree $k$.  It is not hard to see that the homotopy type of $\Gamma_k(M)$ is independent of $k$.  

\begin{theorem}\label{config}\cite{mcduff}  There are  gluing maps $\gamma_k : C_k(M) \to C_{k+1}(M)$ and a  family of maps $\alpha_k : C_k(M) \to \Gamma_k(M)$ that satisfy the following properties:
\begin{enumerate}
\item The induced map in homology $(\alpha_k)_* : H_q(C_k(M)) \to H_q(\Gamma_k(M))$ is an isomorphism if $k$ is sufficiently large. 
\item
These homomorphisms are compatible in the sense that the  following diagrams commute:
$$
\begin{CD}
H_q(C_k(M)) @>(\alpha_k)_* >> H_q(\Gamma_k(M)) \\
@V\gamma_k VV  @VV\cong V \\
H_q(C_{k+1}(M)) @>(\alpha_{k+1})_* >> H_q(\Gamma_{k+1}(M)) 
\end{CD}
$$
where the right vertical map is induced by a homotopy equivalence $\Gamma_k(M)  \xr{\simeq} \Gamma_{k+1}(M)$. 
\end{enumerate}
\end{theorem}

Observe that this theorem answers both Stability Questions 1 and 2 for these moduli spaces.  Question 1 is answered because this theorem  says that the gluing maps $\gamma_k : C_k(M) \to C_{k+1}(M)$ induce isomorphisms in homology through a range of dimensions.  In fact it is proved that the maps $\gamma_k$ induce monomorphisms in homology in \sl all \rm dimensions. Question 2 was answered because this theorem implies  the following.

 Let $C(M)$ be the (homotopy) colimit, $C(M)  = \lim_{k\to \infty} C_k(M)$, where the limit is taken with respect to the gluing maps $\gamma_k$.  
 
 \begin{corollary} The maps $\alpha_k$ induce a map
 $$
 \alpha : \bz \times C(M) \to  \Gamma (M)
 $$
 which induces an isomorphism in homology.
 \end{corollary}

 Notice that if $M$ has a trivial tangent bundle, then $\Gamma (M) \cong C^\infty_{cpt}(M, S^n)$, the space of smooth maps with compact support.  This in turn 
 is homotopy equivalent to $Map_\blt ((M \cup \infty),  S^n)$, where $Map_{\blt}$ denotes the space of basepoint preserving continuous maps.  Moreover, if $\bar M$ is a compact manifold with boundary, having $M$ as its interior, this space can be viewed as the space of maps of pairs,  $Map((\bar M, \p \bar M),  (S^n, \infty))$.

 \med
 An important special case of this theorem, which was proved prior to the proof of Theorem \ref{config} is when $M = \br^n$.  One then has the following well known theorem  about configurations of points in Euclidean space \cite{segal} \cite{may}.
 \begin{theorem}\label{rn}  There are maps $\alpha_k : C_k(\br^n) \to \Omega^n_kS^n$ with the following homological properties:
 \begin{enumerate}
 \item  $(\alpha_k)_* : H_*( C_k(\br^n)) \to H_*(\Omega^n_kS^n)$ is a monomorphism in all dimensions.
 \item $(\alpha_k)_* : H_q( C_k(\br^n)) \to H_q(\Omega^n_kS^n)$ is an isomorphism is $k$ is sufficiently large with respect to $q$.
 \item $\alpha : \bz \times C(\br^n)  \to \Omega^nS^n$  induces an isomorphism in homology.
 \end{enumerate}
 Here $\Omega^n_kS^n$ is the space of (basepoint preserving) self maps of $S^n$ of degree $k$.  
 \end{theorem}

 \subsection{Symmetric groups and braid groups}
  Two special cases of Theorem \ref{rn}  are worth pointing out.  First we consider the case when $n=\infty$.  In this case the space of \sl ordered \rm configurations of points, $F_k(\br^\infty )$ are contractible.  To see this, one considers the projection fibrations, $p_k: F_k(\br^m) \to F_{k-1}(\br^m)$ given by projecting onto the first $k-1$ coordinates.  The fiber of this fibration is $\br^m - \{k-1\}$,  Euclidean space with $(k-1)$-points removed. This space has the homotopy type of a wedge of $(k-1)$ spheres of dimension $m-1$, and therefore its homotopy groups are zero  through dimension $m-2$.  An inductive argument (on $k$), then shows that $\pi_q(F_k(\br^m)) = 0$ for $q \leq m-2$.  
  
 We therefore have that $F_k(\br^\infty)$ is contractible, and has a free action of the symmetric group $\Sigma_k$.  Thus $C_k(\br^\infty)$ is a model for the classifying space $B\Sigma_k$.  Therefore its (co)homology is the (co)homology of the symmetric group $\Sigma_k$.  An alternative viewpoint  is that $C_k(\br^\infty)$ is the moduli
 space $B\Diff (M)$,  as considered in the last section, where $M$ is the zero dimensional manifold consisting of $k$- points.    In any case, Theorem \ref{rn}, applied to the case $n = \infty$ gives the following theorem known as the ``Barratt-Priddy-Quillen theorem" \cite{barrattpriddy}.  Notice that it addresses both Stability Questions 1 and 2 in for the moduli space of points.  This result began a fundamentally important line of research regarding the relationship of finite group theory to stable homotopy theory.  This line of research remains quite active today,  more than 35 years after the proof of the Barratt-Priddy-Quillen theorem.  
  
 \med
 \begin{theorem}\label{bqp}  There are maps $\alpha_k : B\Sigma_k \to \Omega^\infty_k S^\infty$  with the following homological properties:
 \begin{enumerate}
 \item  $(\alpha_k)_* : H_*(B\Sigma_k) \to H_*(\Omega^\infty_kS^\infty)$ is a monomorphism in all dimensions.
 \item $(\alpha_k)_* : H_q(B\Sigma_k) \to H_q(\Omega^\infty_kS^\infty)$ is an isomorphism if $k$ is sufficiently large with respect to $q$.
 \item $\alpha : \bz \times B\Sigma_\infty  \to \Omega^\infty S^\infty$  induces an isomorphism in homology.
 \end{enumerate}
 Here $\Omega^\infty_kS^\infty  = \lim_{n\to \infty} \Omega^n_kS^n$.  
 \end{theorem}
 
 Another important special case of Theorem \ref{rn} is when $n=2$.  We note that $F_2(\br^2)$ has the homotopy type of $S^1$,  whose homotopy groups are $\bz$ in dimension one,  and zero in all other dimensions. In other words, $S^1$ is an Eilenberg-MacLane space, $K(\bz, 1)$.    The fibration  $p_k : F_k(\br^2) \to F_{k-1}(\br^2)$ has fiber $\br^2  - \{k-1\}$, which has the homotopy type of a wedge of $k-1$ circles, which is also a $K(\pi, 1)$.  Now an easy inductive argument (on $k$) implies
 that each $F_k(\br^2)$ is a $K(\pi, 1)$, as is the quotient, $C_k(\br^2)$,  for appropriate groups $\pi$.  
 
 In the case of $C_k(\br^2)$, its fundamental group is Artin's braid group, $\beta_k$.  One can easily visualize that a one-parameter family of configurations of $k$ unordered  points in the plane can be identified with  is a braid in $\br^3$.  So $C_k(\br^2)$ is the classifying space of Artin's braid group $\beta_k$.  Moreover, the natural inclusion, $C_k(\br^2) \hk C_k(\br^\infty)$ is a map $B\beta_k \to B\Sigma_k$, the homotopy type of which is determined by the homomorphism $\beta_k \to \Sigma_k$ that sends a braid to the resulting permutation of the ends of the strings.   Furthermore, the covering space $\Sigma_k \to F_k(\br^2) \to C_k(\br^2)$ makes it apparent that $F_k(\br^2)$ is the classifying space for the \sl pure \rm braid group, $P\beta_k$,  which is the kernel of $\beta_k \to \Sigma_k$. 
The special case of Theorem \ref{rn} in the case $n=2$, establishes the close connections between Artin's braid groups and self maps of $S^2$:

\med
 \begin{theorem}\label{braid}  There are maps $\alpha_k :B\beta_k \to \Omega^2_k S^2$  with the following homological properties:
 \begin{enumerate}
 \item  $(\alpha_k)_* : H_*(B\beta_k) \to H_*(\Omega^2_kS^2)$ is a monomorphism in all dimensions.
 \item $(\alpha_k)_* : H_q(B\beta_k) \to H_q(\Omega^2_kS^2)$ is an isomorphism is $k$ is sufficiently large with respect to $q$.
 \item $\alpha : \bz \times B\beta_\infty  \to \Omega^2 S^2$  induces an isomorphism in homology.
 \end{enumerate}
 \end{theorem}

 We end this section by recalling that one of the  applications of these configuration spaces is that they may be viewed as homogeneous spaces in the following sense.  Suppose, like above, that $M$ is the interior of a manifold with boundary  $\bar M$. Let $\Diff(\bar M, \p)$  be the group of diffeomorphisms of $\bar M$ that fix the boundary, $\p \bar M$ pointwise.  If $M$ is oriented, we write $\Diff^+(\bar M, \p)$ to denote the subgroup of diffeomorphisms that preserve the orientation.
 
 Notice that $\Diff(\bar M, \p)$ acts transitively on the configuration space $C_k(M)$,  and the isotropy group of a fixed configuration of $k$ points is the subgroup
 $\Diff (\bar M, \{k\}, \p)$ that fix those $k$-points (as a set).  This gives a homeomorphism from  $C_k(M)$ to  a homogeneous space,
\begin{equation}
 C_k(M) \cong  \Diff(\bar M, \p)/\Diff (\bar M, \{k\}, \p).
 \end{equation}  Similarly, if $M$ is oriented, it can be written as the quotient, $$C_k(M) \cong  \Diff^+(\bar M, \p)/\Diff^+ (\bar M, \{k\}, \p).$$
 
 In the case when $M = D^2$, the open, two-dimensional disk,  then a famous theorem of Smale asserts that $\Diff^+(\bar D^2, \p)$ is contractible.
 Thus the quotient $\Diff^+(\bar D^2, \p)/\Diff^+(\bar D^2, \{k\},  \p) \cong C_k(\br^2)$ is the classifyng  of the diffeomorphism group,
 $C_k(\br^2) \simeq B\Diff^+(\bar D^2, \{k\},  \p)$.  Now as we saw above,   $C_k(\br^2)$ is a $K(\pi, 1)$, which implies that the homotopy groups of the diffeomorphism group $\Diff^+(\bar D^2, \{k\},  \p)$ are zero in positive dimensions.  This is equivalent to saying that the subgroup of diffeomorphisms that are isotopic to the identity is weakly contractible.  In particular this says that the discrete group of isotopy classes of diffeomorphisms $\Gamma (\bar D^2, \{k\},  \p)$ is the fundamental group of $C_k(\br^2)$.  In general, the group of isotopy classes of diffeomorphisms of a surface is known as the \sl mapping class group \rm of that surface.   (See Section 5 for a more complete discussion.)   In particular this says that  the braid group  can be viewed as the mapping class group,
 $$
 \beta_k \cong \Gamma (\bar D^2, \{k\},  \p).
 $$
 Thus Theorem \ref{braid} can be interpreted as a stability result for the homology of these mapping class groups.  Stability theorems for mapping class groups of positive genus surfaces will be the main subject of section 5 below. 
 
\section{Holomorphic curves  and gauge theory}

In this section we discuss more modern stability theorems that lie in the intersection of topology and algebraic and differential geometry.  These are stability theorems
regarding moduli spaces of holomorphic maps,  bundles, and Yang-Mills connections.  

\subsection{Holomorphic Curves}
The first stability theorem regarding moduli spaces of holomorphic curves  was due to Segal \cite{segal2}.    Let $Rat_d(\bcp^m)$ be the space of based rational maps  in $\bcp^m$ of degree $d$.  That is, $Rat_d(\bcp^m)$ consists holomorphic maps
$$
\alpha : \bc \bp^1 \to \bc \bp^m
$$
that take  $\infty \subset \bcp^1 = \bc \cup \infty, $ to $[1, 1, \cdots , 1] \in \bcp^n$, and have degree $d$.  This moduli space is topologized as a subspace of the continuous two fold loop space,
$Rat_d(\bcp^m) \subset \Omega^2_d \bcp^m$.   This space can be described as a configuration space of $(d+1)$-tuples of complex polynomials,  
$$
z \to (p_0 (z),  p_1(z), \cdots , p_d(z))
$$
where the $p_i$'s are all monic polynomials of degree $d$ that don't share a common root. By identifying a monic polynomial of degree $d$ with its
$d$ roots,  Segal considered this space of rational functions as a certain configuration space of points in $\bc$ (the configuration of the roots of all the polynomials), which then allowed him to describe gluing maps,
$$
Rat_d(\bcp^m) \to Rat_{d+1}(\bcp^m).
$$
The following theorem answers both Stability Questions 1 and 2 in this setting.

\med
\begin{theorem}\label{rat} (Segal \cite{segal2})   
  Both the gluing maps
$$
Rat_d(\bcp^m) \to Rat_{d+1}(\bcp^m)
$$
and the inclusion maps 
$$
Rat_d(\bcp^m) \hk \Omega^2_d \bcp^m
$$
are homotopy equivalences through dimension $d(2m-1)$.  Furthermore both of these maps induce monomorphisms in homology in all dimensions.
\end{theorem}

\med
 Notice that all  of the path components of $\Omega^2 \bcp^m$ are homotopy equivalent to each other.  This can be seen
 by applying loop multiplications by a map of degree one $\iota$, and by a map of degree $-1$ $j$, 
 $$
 \times \iota : \Omega^2_d \bcp^m \to \Omega^2_{d+1} \bcp^m   \quad \times j :  \Omega^2_{d+1} \bcp^m  \to  \Omega^2_d \bcp^m.
 $$
 Since $\iota$ and $j$ are homotopy inverse to each other, each of these maps is a homotopy equivalence.  Then the above theorem implies that if $Rat_\infty \bcp^m$ is  the (homotopy) colimit
 $\lim_{d \to \infty} Rat_d (\bcp^m)$, then there is a homotopy equivalence
 \begin{equation}\label{infiniterat}
 \bz \times Rat_\infty (\bcp^m) \simeq \Omega^2 \bcp^m.
 \end{equation}

We remark that the homotopy type and especially the homology of $\Omega^2 \bcp^m$ is fairly well understood.  Studying the canonical circle bundle, $S^1 \to S^{2m+1} \to \bcp^m$,  yields, by an elementary homotopy argument, that $\Omega^2 S^{2m+1} \to \Omega^2_1\bcp^m$ is a homotopy equivalence. Said another way,
$$
\bz \times \Omega^2 S^{2m+1} \simeq \Omega^2 \bcp^m.
$$

The topology of $Rat_d (\bcp^m)$ was further studied by Cohen-Cohen-Mann-Milgram in \cite{CCMM}.  The stable homotopy type of these rational function
spaces was completely determined, and in particular their homologies were calculated explicitly. The case of $m=1$ is particularly interesting, considering
the fact that both $Rat_d (\bcp^1)$ and the classifying space of the braid groups $B\beta_q$  give an approximation of the homology type of $\Omega^2S^2$  (compare Theorem \ref{rat} and Theorem \ref{braid}).  In \cite{CCMM} the following was proved.

\begin{theorem}\label{CCMM}  $Rat_d(\bcp^1)$ and $B\beta_{2d}$ have the same stable homotopy type.  In particular they have isomorphic homologies.
\end{theorem}

Analogues of the stability Theorem \ref{rat} for rational functions with values in Grassmannians, or more general homogeneous spaces were proved by Kirwan, Guest, and Gravesen in \cite{kirwan},  \cite{guest}, and \cite{gravesen}.  In \cite{cohenjonessegal},  Cohen, Jones, and Segal gave a Morse theoretic proof of Gravesen's theorem, and studied the general question of when a closed, simply connected, integral symplectic manifold has this type of stability property for its moduli space of (based) rational maps
(i.e. holomorphic maps from $\bcp^1$).   The explicit homology type of these more general rational function spaces were computed by Boyer-Hurtubise-Mann-Milgram in \cite{bhmm}.  

Segal also proved a stability result for spaces of holomorphic maps from a higher genus Riemann surface to $\bcp^m$.  Let $\Sigma_g$ be a closed Riemann surface of genus $g$, and let $Hol_d(\Sigma_g, \bcp^m)$  be the space of based holomorphic maps of genus $g$.  Like in the case of rational functions, this is topologized as a subspace of the space of continuous based maps, $Map_d(\Sigma_g, \bcp^n)$.

\begin{theorem}\label{grat} (Segal \cite{segal2})  If $g>0$ the inclusion
$$
Hol_d(\Sigma_g, \bcp^m)  \hk Map_d(\Sigma_g, \bcp^m)
$$
is a homology equivalence up to dimension $(d-2g)(2m-1)$.  
\end{theorem}
Again, it is easy to see that the homotopy type of $Map_d(\Sigma_g, \bcp^m)$ is independent of $d$, and so this theorem describes the stable homology type of $Hol_d(\Sigma_g, \bcp^m)$. 

\med
Segal's theorem can be extended to involve families of complex structures on the surface $\Sigma_g$.  Namely,  let $\cm_{g,d}(\bcp^m)$ be the moduli space of holomorphic curves of genus $g$ and degree $d$  in $\bcp^m$.  More specifically,  $\cm_g(\bcp^m)$ is  defined as follows.  Fix a smooth,  closed, oriented surface $F_g$ of genus $g$. Then
$\cm_g(\bcp^m)$ is the quotient space
\begin{align}
 \cm_g(\bcp^m) &= \{(J, \phi),\, \text{where $J$ is an (almost) complex structure on $F_g$, and} \notag \\
&\phi : (F_g, J) \to \bcp^n \, \text{is holomorphic of degree $d$}\}/ \Diff^+(F_g). \notag
\end{align}
 Here $\Diff^+(F_g)$ is the space of orientation preserving diffeomorphisms which acts diagonally on the space of (almost) complex structures on $F_g$, and on the space
 of maps $F_g \to \bcp^m$.     One can also define a topological analogue,
  $\cm_{g,d}^{top}(\bcp^m)$ which is defined similarly, except that $\phi : F_g \to \bcp^n$ need only be a continuous map.   Recently, D. Ayala proved the following extension of Segal's theorem:
  
  \begin{theorem}\label{ayala}(Ayala \cite{ayala})   The obvious inclusion
  $$
  \cm_{g,d} (\bcp^m)  \hk \cm_{g,d}^{top}(\bcp^m)
  $$
  induces an isomorphism in homology with coefficients in a field of characteristic zero  in dimensions less than  $(d-2g)(2m-1)$. 
  \end{theorem}
  
  The hypothesis that the coefficient field have characteristic zero has to do with the fact that the action of the diffeomorphism group on the space of complex structures has stabilizer groups which are of the homotopy type of finite groups.  This hypothesis can be removed if one defined these moduli spaces using the \sl homotopy \rm orbit spaces of the diffeomorphism groups, rather than the actual orbit spaces.  Equivalently, one could define these moduli spaces as the \sl quotient stack \rm of this action.

This theorem can be viewed as addressing Stability Question 1 in this setting.   Combined with a theorem of Cohen and Madsen \cite{cohenmadsen}, which gives an explicit calculation of $H_*( \cm_{g,d}^{top}(X))$, for $X$ \sl any \rm simply connected space,  through dimension $(g-5)/2$,  this theorem can also be viewed as addressing Stability Question 2.  The Cohen-Madsen result is closely related to, and uses in its proof,  the  work of Harer and Ivanov on the homological stability of mapping class groups \cite{harer}, \cite{ivanov}, and of Madsen and Weiss on their proof of Mumford's conjecture on the stable cohomology of the moduli space of Riemann surfaces \cite{madsenweiss}. These stability theorems will be discussed in more detail in Section 5. 

\subsection{Gauge theory}
\subsubsection{Flat connections on Riemann surfaces}
In a seminal paper  \cite{atiyahbott},  Atiyah and Bott studied the topology of the moduli spaces of Yang-Mills connections on Riemann surfaces, and related them
to moduli spaces of holomorphic bundles.  We will describe one of their main results, and interpret it as a stability theorem for these moduli spaces.
 
Let $\Sigma$ be a closed Riemann surface of genus $g$, and let $E \to \Sigma$ be a principal $G$-bundle, where $G$ is a compact Lie group.  To make the statements of the following theorems easier, we will assume that $G$ is semisimple.  Let $\frak{g}$ be the Lie algebra, and 
Let $ad (E) = E \times_{G} \frak {g} \to \Sigma$ the corresponding ``adjoint bundle",  where $G$ acts on $\frak{g}$ by conjugation.

  Let $\ca (E)$ be the space of connections on $E$, and let $\ca_F(E)$ be the subspace consisting of flat connections.  In the semisimple setting, these flat connections minimize the Yang-Mill functional, 
\begin{align}\label{ym}
\cy \cm  : \ca (E) &\la \br   \notag \\
\cy \cm (A)  &=    \|F_A \|^2  \notag
\end{align}
where $F_A$ is the curvature $2$-form, and $\| \|$ is the $L^2$-norm on $\Omega^*(\Sigma; ad (E))$.  

Let $\cg (E)$ be the \sl gauge group \rm of the bundle $E$.  This is the group of principal $G$-bundle automorphisms of $E \to \Sigma$ that live over the identity map of $\Sigma$.  The inclusion of flat connections, $\ca_F(E) \hk \ca (E)$ is a $\cg$-equivariant embedding, and the following is one of the main results of \cite{atiyahbott}.

\begin{theorem}\label{atiyahbott} (Atiyah-Bott).  The inclusion
$$\ca_F(E) \hk \ca (E)$$
induces an isomorphism on $\cg$-equivariant homology in dimensions less than $2(g-1)r$, where $g$ is the genus of $\Sigma$, and $r$  is the smallest number of the form $\frac{1}{2}dim(G/Q)$, where $Q \subset G$ is any proper, compact subgroup of maximal rank.
\end{theorem}

\bf Remark. \rm The theorem is stated in \cite{atiyahbott} in a slightly different form.  They observe that a connection on $E$ determines a holomorphic structure on the complexification $E^c$.   This allows for the identification of the space  $\ca (E)$ of connections on $E$ with the space of holomorphic structures on $E^c$.  They then show that there is a Morse-type $\cg$-equivariant  stratification of this space of holomorphic bundles,  and they compute the relative codimensions of the strata.  The space of flat connections  is homotopy equivalent to the stratum of ``semi-stable" holomorphic bundles, and the knowledge of the codimension of this stratum in the next lowest stratum (in the partial ordering) leads to a simple
calculation of the ($\cg$-equivariant) connectivity of the inclusion $\ca_F(E) \hk \ca (E)$.  See \cite{cohengalatiuskitchloo} for details of this calculation.

\med
Notice that this theorem can be viewed as a stability theorem since the range of equivariant homology isomorphism increases linearly with the genus.  This can be seen slightly more explicitly as follows.   Recall that if $X$ is a space with an action of a group $K$,  its equivariant homology, $H^*_K (X)$ is defined to be the (ordinary) homology of its \sl homotopy orbit space, \rm $X //K$ defined to be $EK \times_K X$,  where $EK$ is a contractible space with a free $K$-action.
We therefore may consider the following homotopy orbit spaces of the gauge group action,   
$$
\cm_F(E) = \ca_F(E)// \cg  \quad  \cb (E) = \ca (E)//\cg.
$$
Theorem \ref{atiyahbott} can then be restated as follows.

\med
\begin{corollary}\label{atiyahbott2}   The inclusion
$$\cm_F(E) \hk \cb (E) $$
induces an isomorphism in homology in dimensions less than $2(g-1)r$.  
\end{corollary}

Notice from the above discussion, that this can be viewed as a statement about the homology of the moduli space of semistable holomorphic bundles over $\Sigma$ of the topological type of the complexified bundle, $E^c$.

Now recall from \cite{atiyahbott}  that $\cb (E)$ is homotopy equivalent to the mapping space $Map_E(\Sigma, BG)$,  where $BG = EG/G$ is the classifying space of principal $G$-bundles, and $Map_E$  represents the component of the continuous mapping space consisting of maps that classify bundles isomorphic to $E$.  This mapping space has easily described homotopy type (see \cite{atiyahbott}), so this interpretation of the Atiyah-Bott theorem can be viewed as an answer to Stability Question 2 in this setting.

\med
We remark that the Atiyah-Bott theorem has been extended to allow  the complex structure on $\Sigma$ to vary over moduli space.  This was accomplished in \cite{cohengalatiuskitchloo}.  The moduli space under study in that work was defined to be  

$$
   \cm_{g, E}^G = (\ca_F(E) \times J(\Sigma)) // \mathrm{Aut}(E),
$$
where $\mathrm{Aut}(E)$ is the group of $G$-equivariant maps $E \to E$
which lie over some orientation preserving diffeomorphism $\Sigma \to \Sigma$.  
By forgetting the bundle data there is a fibration sequence
$$
\cm_F(E) \to   \cm_{g, E}^G \xr{p} \cm_g
$$
where $\cm_g$ is the moduli space $J(\Sigma) // \Diff^+(\Sigma)$.

The following was proved in   \cite{cohengalatiuskitchloo}.  

\begin{theorem}\cite{cohengalatiuskitchloo}  
There is a map
$$
\alpha :  \cm_{g, E}^G  \to Map_E (\Sigma, BG) // \Diff^+(\Sigma)
$$
that induces an isomorphism in homology  in dimensions less than $2(g-1)r$. 
\end{theorem}
Furthermore  the homology of $Map (\Sigma, BG) // \Diff^+(\Sigma)$ has been computed   in dimensions less than or equal to $(g-4)/2$  explicitly by Cohen and Madsen in \cite{cohenmadsen}.  In particular its rational cohomology (in this range) is freely generated by $H^*(BG)$, and by the Miller-Morita-Mumford canonical classes $\kappa_i$.  Again, this result makes heavy use of Madsen and Weiss's proof of Mumford's conjecture, which will be discussed further in section 5. 

\med
\subsubsection{Self dual connections on four-manifolds and the Atiyah-Jones Conjecture}
One of the most important gauge theoretic    stability theorems was proved by Boyer, Hurtubise, Mann, and Milgram \cite{bhmm}.  This stability theorem had to do with the moduli spaces of self dual connections on $SU(2)$-bundles over $S^4$, and was a verification of a well known conjecture of Atiyah and Jones  \cite{AJ}.

The setup for this theorem is the following.  Isomorphism classes of principal $SU(2)$-bundles over $S^4$ are classified by their second Chern class, $c_2 \in H^4(S^4) \cong \bz$.  Let $p_k : E_k \to S^4$ be a principal $SU(2)$-bundle with Chern class $k$.  Let $\ca_k$ be the space of connections on $E_k$, and $\ca_{sd}^k$ the subspace of self dual connections.  Here we are giving $S^4$ the usual round metric.  $\ca_{sd}$ forms the space of minima of the Yang-Mills functional,
$\cy \cm : \ca_k \to \br$ defined by $\cy \cm (A) = \|F_A\|^2$, much like in the Riemann surface case. 

Let $\cg_k$ be the \sl based \rm gauge group of $E_k$.  This is the group of bundle automorphisms $g : E_k \to E_k$ living over the identity map of $S^4$, with the property that on the fiber over the basepoint $\infty \in \br^4 \cup \infty = S^4$,  $g : (p_k)^{-1} (\infty) \to (p_k)^{-1} (\infty)$ is the identity.  The gauge group $\cg_k$ acts freely on $\ca_k$, so its orbit space $\cb_k = \ca_k/\cg_k$ is the classifying space of the gauge group.   A straightforward homotopy theoretic argument originally due to Gottlieb \cite{gottlieb} says that there is a homotopy equivalence,  
$$
\cb_k \simeq \Omega^4_k BSU(2)
$$
where the subscript denotes the component of the space of based maps $\gamma : S^4 \to BSU(2)$ with $\gamma^*(c_2) = k \in H^4(S^4)$.  The fact that $\Omega BG \simeq G$ is true for any group $G$, implies that 
 $\Omega^4  BSU(2) \simeq \Omega^3SU(2) = \Omega^3S^3$, and hence $\cb_k \simeq \Omega^3_kS^3$. 
 
  Now let $\cm_k (S^4) =  \casd^k /\cg_k$ be the moduli space of self dual connections on $E_k$.  The inclusion $\casd^k \hk \ca_k$ defines a map
  $$
  \cm_k(S^4) \hk \cb_k \simeq \Omega_k^3S^3
  $$
  which was studied by Atiyah and Jones in \cite{AJ}.

Using solutions of the self-dual equations due to the physicist  't Hooft,  Atiyah and Jones were able to prove the following stability theorem  in \cite{AJ}:
\begin{theorem}\cite{AJ}
The map
$$H_*({\cm}_k )\to  H_*(\Omega^3_kS^3)$$
is {\sl surjective}  for $* < k-2$.  
\end{theorem}

  Atiyah and Jones then made the following conjectures:
\begin{enumerate}
\item   The inclusion ${\cm}_k \subset \Omega^3_kS^3$ is a homology
isomorphism in dimensions $t\leq q(k)$ for some increasing function $q(k)$
with $\lim_{k\rightarrow \infty}(q(k)) = \infty$.
\item  The range of the surjection (isomorphism) $q = q(k)$ 
can be explicitly determined as a function of $k$.
\item  The homology statements can be replaced by homotopy
statements in both conjectures 1 and 2. 
\end{enumerate}

The last and strongest statement became commonly known as the
Atiyah-Jones conjecture.  While it is easy to construct
maps $j_k:\Omega^3_k S^3
\to \Omega^3_{k+1}S^3,$ which are homotopy equivalences,
there   was, at the time, no obvious analogous  map  $g_k : {\cm}_k \to \cm_{k+1}.$  Later, Taubes defined such gluing
maps analytically \cite{taubes2}.   In particular he showed that
  the following diagram homotopy commutes
$$\begin{CD}
\cm_k    @>g_k >>  \cm_{k+1}\\
 @VVV   @VVV\\ 
\Omega^3_k S^3 @>j_k>> \Omega^3_{k+1}S^3.
\end{CD}
 $$
This diagram permits (homotopy) direct limits and hence a stable version
of the Atiyah-Jones conjecture.  This was verified by Taubes in
  \cite{taubes} by analytically studying the indices of the nonminimal critical points of the Yang-Mills functional. 

\bg
\begin{theorem}\label{taubes}\cite{taubes}  
Let ${\cm}_{\infty}$ be the homotopy direct limit of the 
${\cm}_k$'s under the inclusions $g_k$ and let
\, $\theta:{\cm}_{\infty}  \la  \Omega^3_0 S^3$ be the
direct limit of the inclusions ${\cm}_k \subset \Omega^3_kS^3$. 
Then $\theta$ is a homotopy equivalence.   
\end{theorem}

Notice that this can be viewed as an answer to Stability Question 2 in this setting.  The answer to Stability Question 1 was supplied by 
Boyer, Hurtubise, Mann, and Milgram with their proof of the Atiyah-Jones conjecture in \cite{bhmm}.  

We remark that one can ask the analogous types of stability questions  when one studies connections on principal bundles for different Lie groups,
and on different four dimensional manifolds. Taubes proved the analogue of  Theorem \ref{taubes} in this full generality.  The full extent to which the analogue  the Atiyah-Jones conjecture holds   is still an open question.

\section{ General linear groups, Pseudoisotopies, and $K$-theory}

\subsection{The stable topology of general linear groups and algebraic $K$-theory}  
Let $R$ be a discrete  ring and $GL_n(R)$ the rank $n$-general linear group.  Understanding  the cohomology of the group $GL_n(R)$ is important  in algebra, topology, algebraic geometry, and number theory.    One may view 
$GL_n(R)$ as the  subgroup of $GL_{n+1}(R)$ consisting of matrices that have zeros in all entries of the $(n+1)^{st}$ row and $(n+1)^{st}$ column except the $(n+1) \times (n+1)$ entry, which is a $1$.   This inclusion defines a map on classifying spaces, $\iota_n : BGL_n(R)  \to BGL_{n+1}(R)$.   Let $BGL(R)$ be the (homotopy) direct limit of these maps.  Recall that Quillen defined the algebraic $K$-groups, $K_i(R)$, to be the $i^{th}$ homotopy group  
$$
K_i(R) = \pi_i(BGL(R)^+),
$$
where Quillen's plus construction is a very explicit construction that changes the homotopy type, but does not change the homology.    
In this context, the Stability Questions 1 and 2 were answered by Charney \cite{charney} in the case when $R$ is a Dedekind domain,
when she proved the following.

\begin{theorem} \cite{charney}  For $R$ a Dedekind domain, the induced maps
$$
\iota_n : H_i (BGL_n(R)) \to H_i(BGL_{n+1}(R)) 
$$
are  isomorphisms if $4i+5 \leq n$.   

If $R$ is the ring of integers in a number field,
$$
\iota_n : \pi_i (BGL_n(R)^+) \to \pi_i(BGL_{n+1}(R)^+)
$$
is an isomorphism for $4i+1 \leq n$. 
\end{theorem}

Generalizations of these homological stability theorems were found by Dwyer \cite{dwyer} and van der Kallen \cite{vanderkallen}. The theorem was generalized to wider classes of rings, to certain classes of nontrivial coefficients modules, and the stability ranges were improved.

\subsection{Pseudoisotopies, and Waldhausen's algebraic $K$-theory of spaces}
Let $M^n$ be a smooth, compact manifold, perhaps with boundary.  The group of pseudoisotopies $\cp (M)$ is defined to be the diffeomorphism group,
$$
\cp (M) = \Diff (M \times I  ; \p M \times I \cup M \times \{0\}).
$$
This group  naturally acts on $\Diff (M)$ in the following way.  Consider the homomorphism $\cp (M) \to \Diff (M)$ which maps  $H \in\cp (M)$ to     $H_1 \in \Diff (M)$  defined to be the restriction of $H$ to $M \times \{1\}$.  The action of $H$  on $\Diff (M)$ is given by $Hf = f \circ H_1$.  Two diffeomorphisms $f_1$ and $f_2$ are said to be pseudoisotopic if they lie in the same orbit of this group action.  Notice that $f_1$ and $f_2$ are \sl isotopic \rm if they lie in the same path component of $\Diff (M)$.   In a seminal  paper \cite{cerf}, Cerf addressed the question: ``If $f_1$ and $f_2$ are pseudoisotopic, are they isotopic?".     In \cite{cerf} Cerf proved the following:

\med
\begin{theorem} \cite{cerf}  Let $M$ be a simply connected, $C^\infty$, closed, $n$-dimensional manifold with $n \geq 6$.  Then $\cp (M)$ is connected.  Therefore in this setting, pseudoisotopic diffeomorphisms are isotopic. 
\end{theorem}

The topology 
of the space of pseudoisotopies has been of great interest ever since that time.  In particular, Hatcher and Wagoner \cite{hatcherwagoner} showed that $\pi_0(\cp (M))$ is not necessarily trivial if $M$ is  not simply connected, even when $n \geq 6$.   We will not state precisely the result of their calculations of $\pi_0(\cp (M))$ here, but they are related to the algebraic $K$-theory of the group ring  of the fundamental group, $K_*(\bz [\pi_1 (M)])$.  

There is a natural    ``suspension" map,
$$
\sigma : \cp (M) \to \cp (M \times I)
$$
defined by essentially letting $\sigma (H)$ be $H \times id$.  We say ``essentially" because a smoothing  process must be done to deform $H \times id$ so that it satisfies the requisite boundary conditions.  Let 
$$
\bp (M) = lim_{k \to \infty} P(M \times I^k)
$$
where the limit is a homotopy colimit under the maps $\sigma$.  This space of   ``stable pseudoisotopies" is of great interest, because
Waldhausen proved that it is an infinite loop space  that can be studied $K$-theoretically.  In particular he defined the notion of the ``Algebraic $K$-theory of a space",  $A(X)$.  (Here $X$ can be any space - not necessarily a manifold.) This is the algebraic $K$-theory of the ``ring up to homotopy",  $Q((\Omega X)_+)$, where as above, $\Omega X$ is the loop space of $X$, and the construction $Q(Y)$ is as defined in section 1. The set of  path connected components can be identified with the group ring,
$$
\pi_0(Q((\Omega X)_+)) \cong \bz[\pi_1(X)]
$$
and $Q((\Omega X)_+)$ can itself be viewed as a type of group ring in the appropriate category of infinite loop spaces.  In any case,
the following was one of Waldhausen's major theorems about these spaces.

\begin{theorem} (See \cite{waldhausen})  The space $A(X)$ splits as a product of infinite loop spaces,
$$A(X) \simeq Wh (X) \times Q(X_+)$$ where $Wh(X)$ is referred to as the ``Whitehead space" of $X$.  In particular if    $X$ is a manifold,    $Wh (X)$ has as its two-fold loop space, the space of stable pseudoisotopies, 
$$
\Omega^2 Wh(X) \cong \bp (X).
$$
\end{theorem}

Of course it then became very important to understand how the space of \sl stable \rm pseudoisotopies $\bp (M)$, which, by Waldhausen's theorem can be studied $K$-theoretically,  approximates the original \sl unstable \rm group of pseudoisotopies $\cp (M)$.    Igusa's stability theorem \cite{igusa} answered this very important question.  It can be viewed as an answer to Stability Question 1 in this context, and together with Waldhausen's theorem, we also have an answer to Stability Question 2.

\med
\begin{theorem}\cite{igusa} (Igusa).  The suspension map $ \sigma : \cp (M^n) \to \cp (M^n \times I)$ induces an isomorphism in homotopy groups in dimensions $k$ so long as $n > max ( 2k+7,  3k+4).$
\end{theorem}

\section{The moduli space of Riemann surfaces, mapping class groups, and the Mumford conjecture}

Probably the most basic, important moduli spaces occurring in geometry and topology  are  the moduli spaces $\cm_{g,n}$ of genus $g$ Riemann surfaces with $n$ boundary components.  Their topology has been of central  interest since the 1960's, and has had important applications to algebraic geometry,  low dimensional topology,  dynamical systems, conformal field theory and string theory in physics, and most recently, algebraic topology. 

Recently, Madsen and Weiss \cite{madsenweiss} identified the ``stable topology" of these moduli spaces, while proving a generalization
of a famous conjecture of Mumford.  In this section we describe some of the ingredients of their stabilization theorems, as well as a related new theorem of Galatius, about automorphisms of free groups \cite{galatius}.
\subsection{ Mapping class groups, moduli spaces, and Thom spaces}   
 The moduli spaces $\cm_{g,n}$ can be defined as follows.  Let $\Sigma_{g,n}$ be a fixed smooth, compact, oriented surface of genus $g > 1$, and $n \geq 0$ boundary components.  Let $\ch_{g,n}$ be the space of hyperbolic metrics on $\Sigma_{g,n}$ with geodesic boundary,  such that each boundary circle has length one. The  moduli space is then defined to be
 $$
 \cm_{g,n} = \ch_{g,n}/ \Diff^+(\Sigma_{g,n}, \p \Sigma_{g,n}),
 $$
 where as earlier,  $ \Diff^+(\Sigma_{g,n}, \p \Sigma_{g,n})$ consists of orientation preserving diffeomorphisms that are the identity on the boundary.   Recall that Teichm\"uller space, $\ct_{g,n}$, can be obtained by taking the quotient of $\ch_{g,n}$ by the subgroup  $\Diff^+_1(\Sigma_{g,n}, \p \Sigma_{g,n})$ of diffeomorphisms isotopic to the identity.   The quotient group $\dgn / \Diff^+_1(\Sigma_{g,n}, \p \Sigma_{g,n})$
 is the discrete group of isotopy classes of diffeomorphisms, known as the mapping class group $\Gamma_{g,n}$.  This can be viewed as the group of path components, $\Gamma_{g,n} = \pi_0(\dgn)$.  In particular we then have
 $$
 \cm_{g,n} = \ct_{g,n}/\Gamma_{g,n}.
 $$
 When the surface has boundary (i.e. $n>0$), the action of the mapping class group on $\ct_{g,n}$ is free.  Moreover, since $\ct_{g,n}$ is homeomorphic to Euclidean space, $\ct_{g,n} \cong \br^{6g-6+2n}$,  the moduli space is a classifying space for the mapping class group,
 \begin{equation}\label{modmap}
  \cm_{g,n}  \simeq B\Gamma_{g,n}.
\end{equation}
  When the surface is closed, the action of $\Gamma_{g,n}$ on $\ct_{g,n}$ has finite stabilizer groups.  This implies that for $k$ any field of characteristic zero, there is still a homology isomorphism,
  \begin{equation}\label{char0}
  H_*(\cm_{g,0}; k) \cong H_*(B\Gamma_{g,0} ; k).
\end{equation}
  Furthermore, since the subgroup $\Diff^+_1(\Sigma_{g,n}, \p \Sigma_{g,n})$ is contractible,  one also has that the full diffeomorphism group  $\dgn$ has contractible components.  This implies that the projection on its components $\dgn \to \Gamma_{g,n}$ is a homotopy equivalence, and hence there is an equivalence of classifying spaces,
\begin{equation}\label{diffmap}
  B\dgn \simeq B\Gamma_{g,n}.
\end{equation}
Putting this equivalence together with equivalence (\ref{modmap}) we see that for $n>0$,  the moduli space $\cm_{g,n}$ is homotopy equivalent to $B\dgn$, and therefore classifies smooth $\Sigma_{g,n}$-bundles.

We now assume $n \geq 1$, and we consider group homomorphisms,
\begin{align}\label{suspend}
\sigma_{1,0} : \Gamma_{g,n} &\to \Gamma_{g+1, n}, \quad \text{and} \notag \\
\sigma_{0,-1} : \Gamma_{g,n} &\to \Gamma_{g, n-1} 
\end{align}
defined as follows. Pick a fixed boundary circle $c \subset \p \Sigma_{g,n}$   Consider an embedding 
$e_{g,n} : \Sigma_{g,n} \hk \Sigma_{g+1, n}$ that sends all of the $(n-1)$ boundary circles of $\p \Sigma_{g,n}$ \sl other than $c$ \rm diffeomorphically to boundary circles of $\p \Sigma_{g+1,n}$, and so that  
$$
\Sigma_{g+1, n} = \Sigma_{g,n} \cup_c T
$$
where $T$ is a surface of genus one with two boundary circles, $c$ and $c'$.
 In other words, $\Sigma_{g+1,n}$ is obtained from $\Sigma_{g,n}$ by \sl gluing \rm in the surface of genus one, $T$.  Given an isotopy class of diffeomorphism of $\Sigma_{g,n}$,   $\gamma \in \Gamma_{g, n}$,  the  element $\sigma_{1,0}(\gamma) \in \Gamma_{g+1,n}$ is the isotopy class    defined by extending  a diffeomorphism in the class of $\gamma$ to all of $\Sigma_{g+1,n}$ by letting it be the identity on $T \subset  \Sigma_{g+1,n}$.  This defines the homomorphism $\sigma_{1,0} : \Gamma_{g,n} \to \Gamma_{g+1,n}$.

 The map $\sigma_{0, -1}$ is defined similarly.   Namely one chooses an embedding 
 $\kappa_{g,n} :  \Sigma_{g,n} \hk \Sigma_{g, n-1}$ that send all of the $(n-1)$  boundary circles of $\p \Sigma_{g,n}$ \sl other than $c$ \rm diffeomorphically to  the $n-1$ boundary circles of $\p \Sigma_{g ,n-1}$, and so that  
$$
\Sigma_{g, n-1} = \Sigma_{g,n} \cup_c D
$$
where $D$ is diffeomorphic to the disk $D^2$.    In other words,   $\Sigma_{g ,n-1}$ is obtained from $\Sigma_{g,n}$  by ``capping off" the boundary circle $c \in \p\Sigma_{g,n}$ by attaching a disk.  By extending a representative diffeomorphism of a class $\gamma \in \Gamma_{g, n}$ by the identity on $D \subset \Sigma_{g, n-1}$, one obtains a homomorphism
$\sigma_{0,-1} : \Gamma_{g,n}  \to \Gamma_{g, n-1}$.    

Notice that the homomorphisms $\sigma_{1,0}$ and $\sigma_{0,-1}$  depend  on the isotopy classes of the choice of embeddings $e_{g,n}$ and $\kappa_{g,n}$ respectively, but the following famous theorem of Harer \cite{harer} shows that any such choice induces an isomorphism in homology through a range.  This result can be viewed as an answer to Stablity Question 1 in this context:

\med
\begin{theorem}\label{harer}  \cite{harer}\cite{ivanov}  For $g >1$ and $n \geq 1$, the homomorphisms,
$\sigma_{1,0}$ and $\sigma_{0,-1}$ induce isomorphisms in the homology of the classifying spaces,
\begin{align}
\sigma_{1,0} : H_q(B\Gamma_{g,n}, \bz) &\xr{\cong} H_q(B\Gamma_{g+1,n}, \bz) \notag \\
\sigma_{0, -1} : H_q(B\Gamma_{g,n}, \bz) &\xr{\cong} H_q(B\Gamma_{g,n-1}, \bz) \notag 
\end{align}
for $2q < g -2$. 
 \end{theorem}

\med
\bf Remarks.  \rm 1.  Harer's original theorem did not have as large a stability range as described here.  This range is due to Ivanov \cite{ivanov}.  The stability range has been improved even further by Boldsen \cite{boldsen}.  

2.  Notice that this result holds for $n=1$.  It therefore implies that the homology of the mapping class groups for \sl closed \rm surfaces, $H_q(B\Gamma_{g, 0})$, is independent of $g$  so long as $2q < g-2$.  

3.  This theorem was generalized to include certain families of twisted coefficients by Ivanov \cite{ivanov2},   Cohen-Madsen \cite{cohenmadsen}, and with improved stability ranges by Boldsen \cite{boldsen}. 

\med
Combining this theorem with statements (\ref{modmap}) and (\ref{char0}) above, one has  the following corollary.

\begin{corollary} For $g>1$ and $n\geq 1$,  the homology of the moduli space of Riemann surfaces
$$
H_q(\cm_{g,n}; \bz)
$$
is independent of $g$ and $n$ so long as $2q < g-2$.   This result holds for the moduli space of  closed surfaces $\cm_{g,0}$ as well,
if one takes homology with coefficients in a field $k$ of characteristic zero.
\end{corollary}

\med
These results can be viewed as answering Stability Question 1 in the case of the moduli space of curves.  One of the major recent advances of the subject was the answering of Stability Question 2 in this setting by Madsen and Weiss \cite{madsenweiss} when they proved a generalization of a long standing conjecture of Mumford \cite{mumford} about the stable cohomology of moduli space, or equivalently, of the stable cohomology of mapping class groups.

A way of stating Mumford's conjecture is as follows.
Let $B\Gamma_{\infty, n}$ be the mapping telescope  (homotopy colimit) of the maps on classifying spaces
$$
B\Gamma_{g, n} \xr{\sigma_{1,0}} B\Gamma_{g+1,  n} \xr{\sigma_{1,0}} B\Gamma_{g+2,  n} \xr{\sigma_{1,0}} \cdots
$$
By the Harer stability theorem, $H_*(B\Gamma_{\infty,  n}, \bz)$ is independent of the number of boundary components $n$.  Moreover it is isomorphic to the homology of the  ``infinite genus moduli space",   $H_*(\cm_{\infty, n}; \bz)$  for $n \geq 1$,  and   if one takes coefficients in a field $k$ of characteristic zero, this homology is isomorphic to the homology of the closed mapping class group $H_*(\cm_{g, 0} ; k)$ if the genus $g$ is large with respect to the homological degree.     Mumford's conjecture was about the stable cohomology $H^*(\bgg; k)$ where $k$ is a field of characteristic zero:  

\begin{conjecture}\label{mumford}(Mumford)  \cite{mumford}   The stable cohomology of the mapping class groups is a polynomial algebra,
$$
H^*(\bgg; k) \cong k[\kappa_1, \kappa_2, \cdots,  \kappa_i, \cdots ]
$$
where $\kappa_i \in H^{2i}(\bgg; k)$ is the Miller-Morita-Mumford  canonical class.
\end{conjecture}

\med
The Miller-Morita-Mumford classes \cite{mumford}, \cite{morita}, \cite{miller} can be defined in the following way.  As remarked above (\ref{diffmap})
there is an equivalence of classifying spaces, $B\Gamma_{g,n} \simeq B\dgn$, and so these spaces classify surface bundles  whose structure group
is this diffeomorphism group.  In particular the cohomology of these classifying spaces is the algebra of characteristic classes of such bundles.  So let
$$
\Sigma_{g,n} \to E \to B
$$
be a smooth bundle with $\dgn$ as its structure group.  Consider the vertical tangent bundle, $T_{vert} E \to E$.  The fiber at $y\in E$, consists of those tangent vectors in  $T_y E$  that are tangent to the fiber surface at $y$ of $E \to B$.   This bundle is a two dimensional,  oriented vector bundle (recall the structure group $\dgn$ consists of orientation preserving diffeomorphisms).  Let $e \in H^2(E)$ be its Euler class.  Then the $\kappa$-classes are defined by integrating powers of $e$  along
fibers,
$$
\kappa_i = \int_{\mathit{fiber}}e^{i+1}  \in H^{2i}(B).
$$
 Alternatively, this is the pushforward in cohomology, $\kappa_i = p_{!} (e^{i+1})$.   Because of the naturality of the pushforward (integration) construction, these
 classes define characteristic classes, and therefore lie in $H^*(B\dgn)$.  Since their construction did not depend on the genus of the surface, they actually define stable cohomology classes in  $H^*(\bgg)$.    We remark that these classes can be constructed directly in  $H^*(\cm_{g,n})$, by integrating along fibers as above, in  the canonical $\Sigma_{g,n}$-bundle,
 $$
 \Sigma_{g,n} \to \cm_{g,n}^1 \to \cm_{g,n}
 $$
 where $\cm_{g,n}^1$ is the moduli space of curves with one marked point. 
 
 We remark that Miller proved in \cite{miller} that the induced map $k[\kappa_1, \kappa_2, \cdots,  \kappa_i, \cdots ] \to H^*(\bgg; k)$ is injective. 
 However to prove Mumford's conjecture (that it is an isomorphism),  Madsen and Weiss employed methods of homotopy theory as well as differential topology.
 In \cite{tillmann}, Tillmann considered the Quillen plus construction, applied to the mapping class groups, $B(\Gamma_{g,n})^+$, and their stabilization,
 $\bgg^+$.  As mentioned earlier, this construction does not alter the homology, so understanding the homotopy type of $\bgg^+$ would yield an understanding of the stable cohomology of mapping class groups.  
 
 In \cite{tillmann},  Tillmann proved that $\bgg^+$ is an infinite loop space.  This result is similar in spirit to     Quillen's result that $BGL(R)^+$ is an infinite loop space (used to define higher algebraic $K$-theory).  In homotopy theory, infinite loop spaces define generalized cohomology theories, but it wasn't clear what generalized cohomology theory $\bgg^+$ defined.  Using a homotopy theoretic model of integrating along fibers that stems from Pontrjagin and Thom's famous work on cobordism theory,   Madsen conjectured what the cohomology theory was.  This conjecture
 was studied by Madsen and Tillmann in \cite{madsentillmann},  and was eventually proved by Madsen and Weiss in \cite{madsenweiss}.  Once this cohomology theory was identified, Mumford's conjecture was an immediate consequence,  as was a description, in principle, of  the stable cohomology, $H^*(\bgg; \bz)$ with \sl integer \rm coefficients.  The cohomology $H^*(\bgg; \bz/p)$ was later computed \sl explicitly \rm by Galatius in \cite{galatius1}.  The answer is quite complicated, but it is entirely defined in terms of rather standard objects in homotopy theory (``Dyer-Lashof operations").  
 
 \med
A basic ingredient in the Madsen-Weiss proof is the use of the Pontrjagin-Thom construction to give a homotopy theoretic model for ``integrating along fibers".
In this particular setting, Madsen and Weiss use a different model of the classifying space   $B\Diff^+(\Sigma_g)$.  (Here I am considering closed surfaces $\Sigma_g$,
but there is an analogous construction for surfaces with boundary, that is equally treated in \cite{madsenweiss}.)   As described in the discussion on Whitney's embedding theorem (\ref{whitney}),  the space $Emb(\Sigma_g, \br^\infty)$, is a contractible space, with a free action of $\Diff^+(\Sigma_g)$.   We then have
that the quotient space, $ Emb(\Sigma_g, \br^\infty)/\Diff^+(\Sigma_g) \simeq B\dsg$.  This is the moduli space of subsurfaces of $\br^\infty$  of genus $g$, which we denote by $\cs_g(\br^\infty)$.  

Consider the space of subsurfaces of $\br^N$, $\cs_g(\br^N) = Emb(\Sigma_g, \br^N)/\dsg$.  There is an obvious $\Sigma_g$-bundle over $\cs_g(\br^N)$,
$$
\Sigma_g \to \cs^1_g (\br^N) \xr{p} \cs_g(\br^N),
$$
where $\cs^1_g(\br^N)$ is the space of subsurfaces of $\br^N$ with a marked point.  Now consider the map
$$
p \times \iota :  \cs^1_g (\br^N) \to \cs_g(\br^N) \times \br^N
$$
where $\iota (S, x) = x \in \br^N$.  This is an embedding.  By identifying $\br^N$ with $B_R(0)$,   the ball of radius  $R$ around the origin, we can consider an induced embedding $ \cs^1_g (\br^N) \hk  \cs_g(\br^N) \times B_R(0) \subset  \cs_g(\br^N) \times \br^N$.  Let $\eta$ be a tubular neighborhood.  More specifically this embedding has a normal bundle, $\nu_N$,  which over each surface $(S,x) \in   \cs^1_g (\br^N)$, is the orthogonal complement of the tangent space $T_x S$ in $\br^N$.  One can then extend this embedding to an embedding of an $\eps$-neighborhood of the zero section of $\nu_N$, for sufficiently small $\eps >0$. The image of this embedding is the tubular neighborhood $\eta$.  

 One then has a ``Pontrjagin-Thom collapse" map
$$
\tau:  \cs_g(\br^N) \times \br^N/ ( \cs_g(\br^N) \times (\br^N - B_R(0)))     \la  \cs_g(\br^N) \times \br^N/(( \cs_g(\br^N) \times \br^N)- \eta).
$$
The left hand side can be identified with the $N$-fold suspension, $\Sigma^N (\cs_g (\br^N)_+)$, and the right hand side can be identified   with the Thom space of the normal bundle $\nu_N$.    Notice that this is an oriented, $(N-2)$-dimensional vector bundle.  An easy bundle theoretic exercise shows that the Whitney sum, 
\begin{equation}\label{normal}
\nu_N \oplus T_{vert} \cs_g^1(\br^N) \cong   \cs_g^1(\br^N) \times \br^N
\end{equation}  viewed as the $N$-dimensional trivial bundle.  We can therefore think of $\nu_N$ as the ``vertical normal bundle"  of the projection map $p : \cs^1_g (\br^N) \to  \cs_g(\br^N)$.  Furthermore this isomorphism and the orientation of $T_{vert} \cs_g^1(\br^N)$ induces an orientation on $\nu_N$.

  The Pontrjagin-Thom map can then be viewed as a map
\begin{equation}\label{thom}
\tau : \Sigma^N (\cs_g (\br^N)_+)  \to Thom (\nu_N),
\end{equation} where we are using the notation $Thom (\zeta)$ to denote the Thom space of a vector bundle $\zeta$.  

Since $\nu_N$ is oriented  there is a Thom isomorphism,  $H^q( \cs_g^1(\br^N)) \xr{\cong} H^{q+N-2} (Thom (\nu_N))$, as well as a suspension isomorphism,
$H^j (\cs_g (\br^N)) \xr{\cong} H^{j+N}( \Sigma^N (\cs_g (\br^N)_+))$.  With respect to these isomorphisms, the induced cohomology homomorphism defined by the  Pontrjagin-Thom map, 
$$
\tau^*: H^*(Thom (\nu_N)) \to H^*(\Sigma^N (\cs_g (\br^N)_+))
$$
induces a homomorphism,
$$
H^q( \cs_g^1(\br^N))  \to H^{q-2}(\cs_g (\br^N))
$$
which is well known to be equal (up to sign) to the fiberwise integration map (or pushforward map) $p_! : H^q( \cs_g^1(\br^N))  \to H^{q-2}(\cs_g (\br^N))$.

This homotopy theoretic view of fiberwise integration has, in some sense a universal model.  Namely, if $Gr^+_2 (\br^N)$ is the Grassmannian of oriented $2$-dimensional subspaces of $\br^N$,   notice that there is a natural map
\begin{align}
j :  \cs^1_g (\br^N)  &\to Gr^+_2 (\br^N) \notag \\
(S, x) &\to T_{vert, x}  \cs^1_g (\br^N) \notag
\end{align}
Notice that the vertical tangent space at $x$   is a subspace of the tangent space, which in turn is a subspace of $\br^N$ since $S \subset \br^N$.  Notice furthermore, that by definition,  if $\gamma_{2,N} \to Gr^+_2 (\br^N)$ is the canonical, oriented $2$-dimensional bundle, then  $j^*(\gamma_{2,N}) = T_{vert} \cs^1_g (\br^N).$  (Recall $\gamma_{2,N}$ consists of pairs, $(V, v)$, where $V \subset \br^N$ is an oriented, $2$-dimensional subspace, and $v \in V$.)  Let $\gamma_{2,N}^\perp$ be the orthogonal complement bundle.  This is the $(N-2)$-dimensional bundle over $Gr^+_2 (\br^N)$ that consists of pairs $(V, w)$, where $V \subset \br^N$ is an oriented, $2$-dimensional subspace, and $w \in V^\perp$.  The  bundle equation (\ref{normal}) induces an isomorphism,  $$j^*(\gamma_{2,N}^\perp) \cong \nu_N.$$
Furthermore $j$ induces a map of Thom spaces,
$$
j: Thom (\nu_N) \to Thom (\gamma_{2,N}^\perp).
$$
The adjoint of the Pontrjagin-Thom map $\tau : \Sigma^N (\cs_g (\br^N)_+)  \to Thom (\nu_N)$, is a map
$\tau : \cs_g(\br^N) \to \Omega^N (Thom (\nu_N))$, and if we compose with the map $j$, we obtain a map
$$
\alpha_{g,N} :  \cs_g(\br^N)  \to \Omega^N(Thom (\gamma_{2,N})^\perp).
$$

Now as observed in \cite{madsenweiss} there are natural inclusions $\Omega^N(Thom (\gamma_{2,N})^\perp)   \hk \Omega^{N+1}(Thom (\gamma_{2,N+1})^\perp)$
that are compatible with the inclusions $ \cs_g(\br^N) \hk  \cs_g(\br^{N+1})$.   We write $\Omega^\infty (Thom (-\gamma_2))$ as the (homotopy) direct limit of these maps.  In the language of homotopy theory, this is the zero space of the Thom spectrum of the virtual bundle $-\gamma_2$,  where $\gamma_2 \to  Gr^+_2 (\br^\infty)$ is the canonical oriented $2$-dimensional bundle.    Notice that this can be identified with the canonical complex line bundle $L \to \bc \bp^\infty$,  and so $\Omega^\infty (Thom (-\gamma_2))$ can be identified with $\Omega^\infty (Thom (-L))$.  (Besides the notation given here, there are several   ``standard" notations for this infinite loop space, including $\Omega^\infty ((\bcp^\infty)^{-L})$, $\Omega^\infty (\bcp^\infty_{-1})$, and more recently,  $\Omega^\infty MTSO(2)$.) 

 In any case, by passing to the limit one has a map
 
\begin{equation}\label{alphag}
 \alpha_g:  B\Diff^+(\Sigma_g) \simeq  \cs_g(\br^\infty)  \to  \Omega^\infty Thom (-\gamma_2).
\end{equation}

The following is the Madsen-Weiss theorem, which supplies a dramatic answer to Stability Question 2 in this setting.

\begin{theorem}\label{madsenweiss}\cite{madsenweiss}   The maps $\alpha_g$ defined above extend to a map
$$
\alpha : \bz \times B\Gamma_{\infty, 1}^+ \to \Omega^\infty Thom (-\gamma_2)
$$
which is a homotopy equivalence (of infinite loop spaces).  In particular, the stable cohomology of the mapping class groups,
$H^*(\bz \times B\Gamma_{\infty, 1}; G)$ is isomorphic to $H^*(\Omega^\infty Thom (-\gamma_2); G)$ for any coefficient group $G$.
\end{theorem}

\med
As mentioned above, the homotopy type of $ \Omega^\infty Thom (-\gamma_2)$ is rather complicated, but it is a natural homotopy theoretic construction, whose basic ingredient is the canonical line bundle $L \to \bcp^\infty$.  In particular the rational cohomology calculation is rather easy, and it is easily seen to imply Mumford's conjecture (\ref{mumford}).  An important implication of this homotopy equivalence and Galatius's calculation \cite{galatius}, is that the stable cohomology of the mapping class groups (and the moduli spaces of curves) has a rich torsion component that classical geometric techniques did not detect. 

\med
Aside from the Pontrjagin-Thom construction, the main idea in the Madsen-Weiss proof was to give a geometric interpretation of the statement in the above theorem. This was done by comparing  concordance (cobordism) classes
of surface bundles,  $M^{n+2} \to X^n$, which are classified by $B\Diff^+(F_g)$  for some $g$, with concordance classes of smooth proper maps $q : M^{n+2} \to X^n$
that come equipped with bundle epimorphisms, $\delta q : TM \times \br^i  \to q^*(TX) \times \br^i$ that live over $q : M \to X$.  Notice that no assumption is made that
the bundle map $\delta q$ is related to the differential $dq$.   Pontrjagin-Thom theory says that as $i$ gets large,  this latter set of concordance classes of maps is classified by
$ \Omega^\infty Thom (-\gamma_2)$.  The comparison of these two sets of concordance classes of maps was studied using an ``$h$-principle" proved by Vassiliev \cite{vassiliev}.   A more detailed outline of the methods used by Madsen and Weiss is contained in the introduction to their paper \cite{madsenweiss}.

\med
A significant simplification of the proof of the Madsen-Weiss theorem was achieved recently by Galatius-Madsen-Tillmann-Weiss \cite{GMTW}.  This paper is about the topology of  ``cobordism categories".   An $n$-dimensional cobordism category has objects consisting of closed $(n-1)$-manifolds, and its morphisms are $n$-dimensional cobordisms between them.  These manifolds may  carry prescribed structure on their tangent bundles, such as orientations,  almost complex structures, or framings.  Much care is given in \cite{GMTW} to give precise definitions to these   topological categories.   Such cobordism categories, aside from the relevance to the stable topology of moduli spaces,  also are relevant in studying topological quantum field theories, and hence their topologies (i.e. the topology of their classifying
spaces) is of great interest.  

In \cite{tillmann}  Tillmann proved that $\bz \times \bgg^+$ has the same homotopy type as the classifying space of  the $2$-dimensional oriented cobordism category $\Omega BCob_2^{or}$.   Her proof involved a clever use of Harer's stability Theorem \ref{harer} \cite{harer}, and a modification of the ``group completion" techniques of McDuff and Segal \cite{mcduffsegal}.
In this remarkable paper, the four authors of \cite{GMTW} then identified the homotopy type of the classifying space of \sl any \rm such cobordism category, in any dimension.  Together with Tillmann's theorem, this gave a simplified proof of the Madsen-Weiss theorem. Moreover,  this theorem, being proved in the generality it was, has had significant influence on the subject beyond the study of the moduli space of curves.  Cobordism theory has been central in differential topology since the original works of Pontrjagin and Thom.  The results of this  paper follow the spirit of Thom's classification of cobordism classes of manifolds, but they go further.  The paper gives
a coherent way of studying the cobordisms that defined the equivalence relation in Thom's theory.  This work has  inspired  considerable work by many people in algebraic and differential topology over the last few years.  Unfortunately, the description of much of this new work is beyond the scope of this paper.
 
\subsection{Automorphisms of free groups}  
One last stability phenomenon that we will discuss concerns automorphisms of the  free group on $n$-generators, $Aut(F_n)$, and the outer automorphism groups, $Out(F_n)$, defined to be the quotient $Out(F_n) = Aut (F_n)/ Inn (F_n)$, where $Inn(F_n) < Aut (F_n)$ is the subgroup of inner automorphisms.   The stability theorems regarding these groups run parallel to, both in statement, and to a certain extent in proof, to  the stability theorems regarding mapping class groups of surfaces due to Harer and Madsen-Weiss, described above.  In particular,  whereas the mapping class group $\Gamma_{g,n}$ is the group of isotopy classes of diffeomorphisms of a surface,  the automorphism group $Aut(F_n)$ is the group of (based) homotopy classes of (based) homotopy equivalences of a graph
$\cg_n$, whose fundamental group is the free group on $n$-generators.  Similarly $Out (F_n)$ can be viewed as the group of unbased homotopy classes of unbased  homotopy equivalences of $\cg_n$.  

In \cite{cullervogtmann} Culler and Vogtmann described a simplicial complex whose simplices are indexed by graphs having fundamental group $F_n$.
This space has a natural action of $Out (F_n)$  and in many ways is analogous to Teichm\"uller space, with its action of the mapping class group.  This
space became known as  ``Outer Space", and has led to many important calculations.  In particular, if one quotients by the action of $Out (F_n)$ one is studying the moduli space of graphs, and it is shown to have the same rational homology as the classifying space, $BOut (F_n)$.  This should be viewed as the analogue of the relationship between the moduli space of curves and the classifying space of the mapping class group (\ref{modmap}) (\ref{char0}). 

Now the natural inclusion $F_n < F_{n+1}$  defines a map $\iota_n : Aut (F_n) \to Aut (F_{n+1})$. We use the same notation for the induced map of classifying spaces, $\iota_n : BAut(F_n) \to BAut(F_{n+1})$.   Similarly, the projection maps $p_n : Aut (F_n) \to Out (F_n)$ define maps on classifying spaces,
$p_n :BAut (F_n) \to BOut (F_n)$.   The following stability theorem of Hatcher and Vogtmann was proved in \cite{hatchervogtmann1.5}, \cite{hatchervogtmann2}, \cite{HVW}. It can be viewed as the analogue of Harer's stability Theorem \ref{harer} above, and can also be viewed as an answer to Stability Question 1 regarding the moduli space of graphs.  

\begin{theorem}\cite{hatchervogtmann1.5},\cite{hatchervogtmann2},\cite{HVW}.  The induced maps in homology,
$$
\iota_* : H_i (BAut (F_n); \bz) \to H_i(BAut (F_{n+1}); \bz)
$$ and
$$
(p_n)_* : H_q(BAut (F_n); \bz) \to H_q (BOut(F_n); \bz)
$$
are isomorphisms for $2i+2 \leq n$ and $2q+4 \leq n$ respectively.
\end{theorem}

This theorem naturally leads to the  problem of computing  stable homology of groups $Aut(F_n)$.   In \cite{hatcher} Hatcher   conjectured
that the rational stable homology is zero.  More precisely,  let $Aut(F_\infty)$ be the direct limit of the groups,  $$Aut (F_\infty) = \lim_{n\to\infty} Aut(F_n).$$

\begin{conjecture} \cite{hatcher}(Hatcher)  The rational homology groups,
$$
H_i(BAut(F_\infty); \bq) = 0
$$
for $i>0$.
\end{conjecture}

This conjecture was recently proved, in dramatic fashion, by S. Galatius in \cite{galatius}. Galatius actually proved a theorem that computes this stable homology with \sl any \rm coefficients.

\med
\begin{theorem}\label{galatius}\cite{galatius}(Galatius)   Let $\Sigma_n$ be the symmetric group on $n$-letters.  View $\Sigma_n$ as the subgroup of $Aut(F_n)$  given by permutations of the generators of $F_n$.  Then the map on classifying spaces $B\Sigma_n \to  BAut (F_n)$
induces an isomorphism in homology,
$$
H_i(B\Sigma_n; G)   \xr{\cong}  H_i(BAut (F_n); G)
$$
for $2i+2 \leq n$, and $G$ any coefficient group.
In particular the induced map
$$
B\Sigma_\infty \to BAut(F_\infty)
$$
is a homology equivalence.     When one applies the Quillen plus construction, there are homotopy equivalences,
$$
\bz \times B\Sigma_\infty^+ \xr{\simeq}  \bz \times BAut (F_\infty)^+ \xr{\simeq}  \Omega^\infty S^\infty
$$
where, like above, $\Omega^\infty S^\infty= \lim_{n\to \infty} \Omega^nS^n.$
\end{theorem}

\med
Notice that, since the symmetric groups are finite, they have trivial rational homology.  Thus  Hatcher's conjecture is a corollary of Galatius's theorem.  Now the homology of the symmetric groups is completely known with any field coefficients, and hence the stable homology of the automorphism groups of free groups is similarly now known.  Notice this result is compatible with the Barratt-Priddy-Quillen Theorem \ref{bqp} regarding the stable homology of symmetric groups.

Galatius's  argument  is similar to the Madsen-Weiss argument in spirit, but involved
many new ideas and constructions.     A key idea in Galatius's argument is to build a model for $BOut(F_n)$ as a space of graphs embedded in Euclidean space.
This builds on the Culler-Vogtmann model of ``Outer Space".  He then defined a ``scanning procedure",  much like what was used by Segal in his study of rational functions \cite{segal2},   to define a map
$$
\alpha : BOut(F_n) \to \Omega^\infty \Phi
$$ 
where $\Omega^\infty \Phi$ is the natural home for the image of a Pontrjagin-Thom type map.  (Notice that this is not the standard Pontrjagin-Thom constuction, since the graphs involved are obviously not smooth manifolds.)
This space is itself defined from a sheaf of (noncompact) graphs.  The induced maps $BAut (F_n) \to \Omega^\infty \Phi$ then extends
to a map $\alpha : BAut (F_\infty)^+ \to \Omega^\infty \Phi$.  This map is analogous to the Madsen-Tillmann  map $\alpha : B\Gamma_{\infty, 1}^+ \to \Omega^\infty Thom (-\gamma_2)$ described in Theorem \ref{madsenweiss} above.  Galatius then proved that this map is a homotopy equivalence.  Finally he proved that 
$\Omega^\infty \Phi$ has the homotopy type of $\Omega^\infty S^\infty$.  

We remark that Galatius's method of proof is quite general, and in particular leads to a further simplification of the proof of the Madsen-Weiss theorem, as well as the theorems of \cite{GMTW} on cobordism categories.  It has also lead to considerable generalizations of these theorems (see \cite{ayala2}\cite{genauer}\cite{oscarsoren}). 

\section{Final Comments}
The stability theorems considered here come in different types, and have a variety of different characteristic features.  However they have all had a significant impact on their field of research, and in some cases that impact has been quite dramatic.  It therefore seems that it would be quite valuable to understand the common features of the classifying spaces and moduli spaces that admit stability theorems,  and to try to understand common features of their proofs.  

For example, many of the stability theorems considered here have to do with classifying spaces of sequences of groups.  They included   braid groups,
  symmetric groups,   general linear groups,   mapping class groups, and automorphisms of free groups.  The common way in which Stability Question 1 was proved in these cases has been considered by Hatcher and Wahl.  In all of these cases, simplicial complexes with the appropriate group actions were found or constructed,
  with certain criteria on the stabilizer subgroups.  Then, typically, a spectral sequence argument was used to inductively prove a stabilty theorem.  Understanding the general properties  of these group actions (i.e. finding axioms) that would imply these stability theorems has lead to the discovery of new such theorems (see for example, \cite{hatcherwahl}, \cite{wahl}).  
  
 As for other general features, notice that some of the above stability theorems that addressed Stability Question 2  involved some variation of the Pontrjagin-Thom construction.  This  was true of the proofs of Theorems \ref{config}, \ref{rat},  \ref{madsenweiss}, \ref{galatius}  described above.  It would certainly be of great value to understand under what conditions the Pontrjagin-Thom construction yields a (homology) equivalence, and therefore a stability theorem.
 
 Some of the stability theorems described above are, in a sense, more analytic in nature.  They concern the moduli space of solutions to a differential equation, such as the Cauchy-Riemann equation in the case of holomorphic curves, or self-duality equations in the case of Yang-Mills moduli spaces.  In these cases, both Stability Questions 1 and 2 involve understanding  the relative topology of the moduli space of solutions inside the entire configuration space (e.g. holomorphic curves inside all smooth curves,  or self-dual connections inside all connections).  In \cite{cohenjonessegal} Cohen, Jones, and Segal discussed sufficient Morse-theoretic conditions
 on when the space of rational maps to a symplectic manifold approximate the topology of all continuous maps of $S^2$ to the manifold.  Their condition involved
 a kind of homogeneity  property.  However it is far from understood, in general, for what type of symplectic manifolds, and for what choices of compatible almost complex structure is there a stability theorem for spaces of pseudo-holomorphic curves.

Stability theorems have been important in both algebraic and differential topology, as well as both algebraic and differential geometry.  Clearly understanding
the conditions under which they do and do not occur is a research goal of real value in all these areas.

 \end{document}